\documentclass[11pt]{conm-p-l}
\usepackage{amscd,amssymb}
\usepackage{xspace}

\newtheorem{theorem}{Theorem}[section]
\newtheorem{proposition}[theorem]{Proposition}
\newtheorem{lemma}[theorem]{Lemma}

\theoremstyle{definition}
\newtheorem{definition}[theorem]{Definition}

\theoremstyle{remark}

\numberwithin{equation}{section}

\DeclareMathOperator{\dg}{dg-}

\DeclareMathOperator{\Ext}{Ext}
\DeclareMathOperator{\Tor}{Tor}
\DeclareMathOperator{\Hom}{Hom}
\DeclareMathOperator{\cok}{cok}
\newcommand{\cat}[1]{\mathcal{#1}}
\newcommand{\ideal}[1]{\mathfrak{#1}}
\newcommand{\mathcolon}{\colon\,}
\newcommand{\Z}{\mathbb{Z}}
\newcommand{\Q}{\mathbb{Q}}
\newcommand{\ulp}{\textup{(}}
\newcommand{\urp}{\textup{)}}
\newcommand{\llp}{left lifting property with respect to\xspace}
\newcommand{\rlp}{right lifting property with respect to\xspace}
\newcommand{\boxprod}{\mathbin\square}
\newcommand{\uc}{\textup{:}}

\newcommand{\Ch}[1]{\text{Ch} (#1)}

  \copyrightinfo{2006}
    {American Mathematical Society}

\begin{document}

\title{Cotorsion pairs and model categories}

\author{Mark Hovey}
\address{Department of Mathematics \\ Wesleyan University
\\ Middletown, CT 06459}
\email{hovey@member.ams.org}

\subjclass[2000]{Primary 55U35; Secondary 14F05, 16D90, 18E30, 18G55}
\date{\today}
\maketitle
\tableofcontents


The purpose of this paper is to describe a connection between model
categories, a structure invented by algebraic topologists that allows
one to introduce the ideas of homotopy theory to situations far removed
from topological spaces, and cotorsion pairs, an algebraic notion that
simultaneously generalizes the notion of projective and injective
objects.  In brief, a model category structure on an abelian category
$\cat{A}$ that respects the abelian structure in a simple way
is equivalent to two compatible complete cotorsion pairs on $\cat{A}$.
This connection enables one to interpret results about cotorsion pairs
in terms of homotopy theory (for example, the flat cover
conjecture~\cite{bican-el-bashir-enochs} can be thought of as the search
for a suitable cofibrant replacement), and vice versa. 

Besides describing this connection, we also indicate some applications.
The stable module category of a finite group $G$ over a field $k$ is a
basic object of study in modular representation theory; it is a
triangulated category because injective and projective $k[G]$-modules
coincide.  Cotorsion pairs can be used to construct two different model
structures on the category of $K[G]$-modules where $K$ is a commutative
Gorenstein ring (such as $\Z $, for example).  The homotopy category of
these two different model structures is the same; it is a triangulated
category that can reasonably be called the stable module category of $G$
over $K$.  This opens up the possibility of ``integral representation
theory'' along the lines of modular representation theory. 

Another application is due to Jim Gillespie.  In algebraic geometry, a
common object of study is the derived category of a scheme, obtained
from chain complexes of quasi-coherent sheaves by inverting maps that
induce isomorphisms on homology.  There is a derived tensor product on
this derived category if the scheme is nice enough, but the construction
used in algebraic geometry seems to the author to be somewhat ad hoc and
difficult to work with.  The essential difficulty is that there are not
enough projective quasi-coherent sheaves in general.  Gillespie has
proved a general theorem about promoting a cotorsion pair on an abelian
category to a model structure on chain complexes over that category.
When applied to quasi-coherent sheaves, it produces a model structure 
compatible with the tensor product of chain complexes of sheaves.  The
existence of the derived tensor product and its expected properties now
follow formally from this model structure.  

This paper is an expanded version of two talks given by the author at
the Summer School on the Interactions between Homotopy Theory and
Algebra at the University of Chicago, July 26 to August 6, 2004.  For
more details, the reader can consult the
papers~\cite{hovey-cotorsion},~\cite{gillespie-modules},
and~\cite{gillespie-sheaves}.  The connection between model structures
and cotorsion pairs is also discussed by Beligiannis and Reiten
in~\cite[Chapter~VIII]{beligiannis-reiten}. The author would like to
thank the organizers of the Summer School for inviting him to speak.
The author also thanks the referee for pointing out a subtlety with
hereditary cotorsion pairs that the author and Gillespie had both
missed.  
\setcounter{tocdepth}{3}
\tableofcontents

\section{Cotorsion pairs}

Cotorsion pairs were invented by Luigi Salce~\cite{salce} in the
category of abelian groups, and were rediscovered by Ed Enochs and
coauthors in the 1990's. A \textbf{cotorsion pair} in an abelian
category $\cat{A}$ is a pair $(\cat{D},\cat{E})$ of classes of objects
of $\cat{A}$ each of which is the orthogonal complement of the other
with respect to the $\Ext$ functor.  That is, we have 
\begin{enumerate}
\item $D\in \cat{D}$ if and only if $\Ext^{1}(D,E)=0$ for all $E\in
\cat{E}$; and 
\item $E\in \cat{E}$ if and only if $\Ext^{1}(D,E)=0$ for all $D\in
\cat{D}$.  
\end{enumerate}
Cotorsion pairs have been used to study covers and
envelopes~\cite{enochs-lopez-book}, \cite{enochs-jenda-book},
particularly in the proof of the flat cover
conjecture~\cite{bican-el-bashir-enochs}.  They have also been used in
tilting theory~\cite{bazzoni-eklof-trlifaj} and in the representation
theory of Artin algebras~\cite{krause-solberg}.  

The most obvious example of a cotorsion pair is when
$\cat{D}=\cat{A}$, in which case $\cat{E}$ is the class of injective
objects.  Similarly, we could let $\cat{E}=\cat{A}$, in which case
$\cat{D}$ is the class of projective objects.  

Based on this example, we say that a cotorsion pair
$(\cat{D},\cat{E})$ \textbf{has enough projectives} if for all $X$ in
our abelian category $\cat{A}$ there is a short exact sequence 
\[
0 \xrightarrow{} E \xrightarrow{} D \xrightarrow{} X \xrightarrow{} 0
\]
where $D\in \cat{D}$ and $E\in \cat{E}$.  So $\cat{A}$ has enough
projectives in the usual sense if and only if the cotorsion pair
(projectives, everything) has enough projectives.  On the other hand,
the cotorsion pair (everything, injectives) always has enough
projectives.  Dually, we say that $(\cat{D}, \cat{E})$ \textbf{has
enough injectives} if for all $X$ in $\cat{A}$ there is a short exact
sequence
\[
0 \xrightarrow{} X \xrightarrow{} E \xrightarrow{} D \xrightarrow{} 0
\]
with $E\in \cat{E}$ and $D\in \cat{D}$.  If $(\cat{D},\cat{E})$ has
enough projectives and enough injectives, we say that it is a
\textbf{complete cotorsion pair}.  

Perhaps the most useful cotorsion pair, and the one that gives the
subject its name, is the flat cotorsion pair.
Here $\cat{A}$ is the category of $R$-modules for some ring $R$,
$\cat{D}$ is the category of flat $R$-modules, and $\cat{E}$ is what it
has to be, the collection of all modules $E$ such that
$\Ext^{1}_{R}(D,E)=0$ for all flat $D$.  Such modules are called
cotorsion modules.  

It is not at all obvious that this is a cotorsion pair, or what
cotorsion modules look like.  A brief digression may be warranted to
describe this important example.  

First of all, a short exact sequence 
\[
0 \xrightarrow{} A \xrightarrow{} B \xrightarrow{} C \xrightarrow{} 0
\]
is called \textbf{pure} if it remains exact upon applying the functor
$M\otimes_{R}(-)$ for any $R$-module $M$.  My favorite reference for
purity and many other algebraic topics is~\cite{lam}; purity is
discussed in Section~4J, where it is proved, among other things, that
the pure exact sequences are the colimits of split exact sequences, and
that any short exact sequence where the right-hand entry $C$ is flat is
automatically pure.  Purity is of considerable interest to logicians
interested in the model theory of modules~\cite{herzog}.  

A module $A$ is \textbf{pure injective} if every pure exact sequence
with $A$ as the left-hand entry is in fact split.  There are lots of
these around; most importantly, if $M$ is any $R$-module, then
$M^{+}=\Hom_{\Z}(M,\Q /\Z)$ is always pure injective.  (Note that if $M$
is a left $R$-module, then $M^{+}$ is a right $R$-module).  And every
pure injective module is cotorsion (because any short exact sequence
that ends in a flat is automatically pure), so this gives us a source of
cotorsion modules.  Using these facts, we can prove that (flat,
cotorsion) is in fact a cotorsion pair.  Indeed, it suffices to show
that if $\Ext^{1}(D,E)=0$ for all cotorsion $E$, then $D$ is flat.  But
this means that $\Ext^{1}(D,M^{+})=0$ for all $M$.  Using the derived
version of the $\Hom$ and tensor adjointness, and the fact that $\Q /\Z$
is injective as an abelian group, we see that $\Tor^{1}(D,M)^{+}=0$ for
all $M$, which implies that $D$ is flat.  

It was an open question for a long time whether the cotorsion pair
(flat, cotorsion) was complete.  This became known as the \textbf{flat
cover conjecture}.  It was eventually proved when Eklof and
Trlifaj~\cite{eklof-trlifaj}, working from a mathematical logic
perspective, reinvented and applied some model category theoretic
techniques (though I don't believe they were aware of the connection).
Bican, El Bashir, and Enochs then used the Eklof-Trlifaj result to prove
the flat cover conjecture~\cite{bican-el-bashir-enochs}.

\section{Relation between cotorsion pairs and model categories}

Recall that a \textbf{model category} is a category $\cat{M}$, which I
will assume has all limits and colimits, together with three
subcategories (the \textbf{model structure}) called the \textbf{weak
equivalences}, \textbf{cofibrations}, and \textbf{fibrations} that must
satisfy various axioms.  Model categories allow one to export the
methods of algebraic topology from topological spaces to more general
situations.  For the author, the guiding principle is that anytime one
has a class of maps that are not isomorphisms but that one wishes were
isomorphisms, there should be a model category lurking in the background
for which those maps are the weak equivalences.  Furthermore, making
that model structure explicit often gives rise to additional structures
that were not readily apparent beforehand.  One of the simplest
interesting examples is the category of (unbounded) chain complexes of
modules over a ring, where the weak equivalences are the homology
isomorphisms.  A model category has a homotopy category, obtained by
formally inverting the weak equivalences; in the case of chain complexes
this homotopy category is known as the derived category of the ring and
is of central importance in homological algebra and algebraic geometry.
When the ring is commutative, one would like a well-behaved derived
tensor product on the derived category.  Trying to construct this
derived tensor product without a model structure can be quite painful,
but with a particular model structure on chain complexes, known as the
projective model structure, the existence and the expected properties of
the derived tensor product follow easily.

For this paper, we will not need more than the introduction to model
categories in~[Sections~1-3]\cite{goerss-model} in these proceedings,
whose notation we will follow.  Another good introduction
is~\cite{dwyer-spalinski}, and reference books on model categories
include~\cite{goerss-jardine}, ~\cite{hirschhorn},
and~\cite{hovey-model}, as well as the original
source~\cite{quillen-htpy}.

Suppose we have a cotorsion pair $(\cat{D},\cat{E})$.  This means that
given any short exact sequence 
\[
0\xrightarrow{}A\xrightarrow{i}B\xrightarrow{}D\xrightarrow{}0
\]
with $D\in \cat{D}$, and for any $E\in \cat{E}$, the map
$\cat{A}(B,E)\xrightarrow{}\cat{A}(A,E)$ is surjective, because
$\Ext^{1}(D,E)=0$.  In model category language, this says that given the
commutative diagram below
\[
\begin{CD}
A @>f>> E \\
@ViVV @VVV \\
B @>>> 0
\end{CD}
\]
we can find a lift $g\mathcolon B\xrightarrow{}E$ making both triangles
commute.  This looks like the lifting axiom for model
categories~\cite[Definition~1.2]{goerss-model}.  That is, if we imagine
$i$ to be an acyclic cofibration (both a cofibration and a weak
equivalence), then $E\xrightarrow{}0$ looks like a fibration, so that
$\cat{E}$ consists of fibrant objects.  This suggests that there should
be some relation between model categories and cotorsion pairs.

\subsection{Abelian model categories}

For this relationship between model categories and cotorsion pairs to
hold, we need some relation between the model structure on $\cat{A}$ and
the abelian structure.

\begin{definition}\label{defn-abelian}
An \textbf{abelian model category} is an complete and cocomplete abelian
category $\cat{A}$ equipped with a model structure such that
\begin{enumerate}
\item A map is a cofibration if and only if it is a monomorphism with
cofibrant cokernel.  
\item A map is a fibration if and only if it is an epimorphism with
fibrant kernel.
\end{enumerate}
\end{definition}

Here we are using the definition of model structure
from~\cite[Section~1.2]{goerss-model}; in particular, we will not assume
the factorizations in the factorization axiom M.5 are functorial, as was
done in~\cite{hovey-model}.  Most of the standard model structures on
abelian categories are abelian model structures.  For example, in the
projective model structure on chain complexes, the cofibrations are the
monomorphisms with cofibrant (=DG-projective) cokernel, the fibrations
are the epimorphisms, and the weak equivalences are the homology
isomorphisms. A complex is DG-projective if each entry is projective,
and if any map from it to an exact complex is chain homotopic to $0$.

A trivial example of a model structure that is not abelian is the one
where weak equivalences are isomorphisms and all maps are cofibrations
and fibrations (this example is not mentioned in~\cite{goerss-model},
but one can easily check it satisfies the axioms).  A less trivial
example is the absolute model structure on chain
complexes~\cite[Example~3.4]{hovey-christensen-relative}, where the weak
equivalences are chain homotopy equivalences, and everything is
cofibrant and fibrant.  In this model structure, the cofibrations are
the degreewise split monomorphisms and the fibrations are the degreewise
split epimorphisms.  So an epimorphism with fibrant kernel is usually
not a fibration.  However, it is possible to modify the definition of
abelian model category to include this example and many others, using
the idea of a proper class of short exact sequences.  This is the same
thing as an additive subfunctor of the $\Ext$ functor, so there is also
a modified definition of a cotorsion pair using this subfunctor.  In the
case of the absolute model structure, our proper class is the class of
degreewise split sequences.

Suppose $\cat{A}$ is an abelian model category, and $p\mathcolon
X\xrightarrow{}Y$ is an acyclic fibration with kernel $K$.  Then
$K\xrightarrow{}0$ is a pullback of $p$, so is also an acyclic fibration.
We say that $K$ is an \textbf{acyclic fibrant object}, and so an acyclic
fibration is an epimorphism with acyclic fibrant kernel.  In fact, the
converse is true as well; in an abelian model category, every
epimorphism with acyclic fibrant kernel is an acyclic fibration.  The
proof of this converse can be found
in~\cite[Proposition~4.2]{hovey-cotorsion}, but it requires
Proposition~\ref{prop-abelian-cotorsion} below.  Dually, an acyclic
cofibration is easily seen to be a monomorphism with acyclic cofibrant
cokernel; the converse is again less obvious but true.

Now does an abelian model structure have something to do with cotorsion
pairs?  Yes!

\begin{proposition}\label{prop-abelian-cotorsion}
Suppose $\cat{A}$ is an abelian model category.  Let $\cat{C}$ denote
the class of cofibrant objects, $\cat{F}$ the class of fibrant objects,
and $\cat{W}$ the class of acyclic objects \ulp those that are weakly
equivalent to $0$\urp .  Then $(\cat{C}\cap \cat{W},\cat{F})$ and
$(\cat{C},\cat{F}\cap \cat{W})$ are complete cotorsion pairs.  
\end{proposition}

Details of the proof of this proposition can be found
in~\cite{hovey-cotorsion}.  

\begin{proof}[Sketch of proof]
We just do the $(\cat{C},\cat{F}\cap \cat{W})$ case as the other is
similar.  There are 5 steps to the argument.
\begin{enumerate}
\item $\Ext^{1}(C,K)=0$ for cofibrant $C$ and acyclic fibrant $K$.
An element of $\Ext^{1}(C,K)$ is represented by a short exact sequence 
\[
0 \xrightarrow{} K \xrightarrow{i} X \xrightarrow{p} C \xrightarrow{} 0.
\]
Since $C$ is cofibrant, $i$ is a cofibration.  By lifting in the diagram 
\[
\begin{CD}
K @= K \\
@ViVV @VVV \\
X @>>> 0
\end{CD}
\]
we get a splitting of our short exact sequence.  

\item If $\Ext^{1}(A,K)=0$ for all acyclic fibrant $K$, then $A$ is
cofibrant.  Prove this by showing that $\cat{A}(A,-)$ takes acyclic
fibrations to surjections, so $0\xrightarrow{}A$ has the \llp acyclic
fibrations (see~\cite[Remark~1.5(2)]{goerss-model} for the definition of
the left lifting property).  

\item If $\Ext^{1}(C,X)=0$ for all cofibrant $C$, then $X$ is acyclic
fibrant.  Prove this by showing $X\xrightarrow{}0$ has the \rlp
cofibrations.  

\item The cotorsion pair has enough projectives.  Prove this by
factoring $0\xrightarrow{}X$ into a cofibration followed by an acyclic
fibration.  
\item The cotorsion pair has enough injectives.  Prove this by
factoring $X\xrightarrow{}0$ into a cofibration followed by an acyclic
fibration.  
\end{enumerate}
\end{proof}

\subsection{From cotorsion pairs to an abelian model category}

We have seen that an abelian model structure gives rise to two
compatible complete cotorsion pairs.  Can we go the other way?  Well,
no, not without some more hypotheses.  Recalling the model category
axioms~\cite[Definition~1.4]{goerss-model}, we have the lifting and
factorization axioms, but we also have the two out of three axiom and
the retract axiom.  To make these other axioms work we are going to need
some hypothesis on $\cat{W}$.

\begin{definition}\label{defn-thick}
A nonempty subcategory of an abelian category is called \textbf{thick}
if it is closed under retracts and whenever two out of three entries in
a short exact sequence are in the thick subcategory, so is the third.
\end{definition}

\begin{lemma}\label{lem-thick}
Suppose $\cat{A}$ is an abelian model category and $\cat{W}$ is the
class of acyclic objects.  Then $\cat{W}$ is thick.  
\end{lemma}

We leave the proof to the reader; it can also be found
in~\cite{hovey-cotorsion}. 

So now we get the desired theorem. 

\begin{theorem}\label{thm-abelian}
Suppose $\cat{C}$, $\cat{F}$, and $\cat{W}$ are three classes of objects
in a bicomplete abelian category $\cat{A}$, such that \begin{enumerate}
\item $\cat{W}$ is thick. 
\item $(\cat{C}, \cat{F}\cap \cat{W})$ and $(\cat{C}\cap
\cat{W},\cat{F})$ are complete cotorsion pairs. 
\end{enumerate}
Then there exists a unique abelian model structure on $\cat{A}$ such
that $\cat{C}$ is the class of cofibrant objects, $\cat{F}$ is the class
of fibrant objects, and $\cat{W}$ is the class of acyclic objects.  
\end{theorem}

The proof of this theorem (which can be found in~\cite{hovey-cotorsion})
is interesting, as it does not follow the usual path for proving
something is a model category.  Usually the main difficulty is proving
the lifting and factorization axioms, but in this case the main
difficulty is defining the weak equivalences and proving the
two-out-of-three axiom, which is usually trivial.

It is clear that we should define $f$ to be a cofibration if $f$ is a
monomorphism with cokernel in $\cat{C}$, a fibration if $f$ is an
epimorphism with kernel in $\cat{F}$, an acyclic cofibration if $f$ is a
monomorphism with cokernel in $\cat{C}\cap \cat{W}$, and an acyclic
fibration if $f$ is an epimorphism with kernel in $\cat{F}\cap \cat{W}$.
But weak equivalences do not have to be monomorphisms or epimorphisms,
so we can't define them in the same way.  Instead, we define $f$ to be a
weak equivalence if it is the composition of an acyclic cofibration
followed by an acyclic fibration.

There are now a great many things to check.  Just to give the flavor of
the argument, we prove a few results needed for
Theorem~\ref{thm-abelian}.

\begin{lemma}\label{lem-composition}
Cofibrations, acyclic cofibrations, fibrations, and acyclic fibrations
are all closed under compositions.  
\end{lemma}

\begin{proof}
Suppose $i\mathcolon A\xrightarrow{}B$ and $j\mathcolon
B\xrightarrow{}C$ are cofibrations.  We have a short exact sequence
\[
0\xrightarrow{}\cok i \xrightarrow{} \cok ji \xrightarrow{}\cok
j\xrightarrow{}0.  
\]
This is a special case of the snake lemma.  Because $\cat{C}$ is the
left half of a cotorsion pair, it is closed under extensions.  Thus
$\cok ji\in \cat{C}$ and so $ji$ is a cofibration.  Because $\cat{W}$ is
thick, if $i$ and $j$ are acyclic cofibrations, so is $ji$.  The
fibration case is similar.
\end{proof}

\begin{proposition}\label{prop-factorization}
Every map $f$ can be factored as $f=qj=pi$, where $j$ is a cofibration,
$q$ is an acyclic fibration, $i$ is an acyclic cofibration, and $p$ is a
fibration.
\end{proposition}

\begin{proof}
This proceeds in stages. The two cases are similar, so we just do the
$qj$ case.  We first assume $f\mathcolon A\xrightarrow{}B$ is a
monomorphism already, with cokernel $C$.  Since $(\cat{C},\cat{F}\cap
\cat{W})$ is a complete cotorsion pair, there is a surjection
$QC\xrightarrow{}C$ where $QC\in \cat{C}$, with kernel $K$ in
$\cat{F}\cap \cat{W}$.  By taking the pullback, we get a monomorphism
$j\mathcolon A\xrightarrow{}B'$ with cokernel $QC$, so $j$ is a
cofibration.  We also get $q\mathcolon B'\xrightarrow{}B$, which is a
surjection with kernel $K$, so an acyclic fibration as required.

Now suppose $f$ is an epimorphism with kernel $K$.  Then we can repeat
the same trick, using an embedding $K\xrightarrow{}RK$ with $RK\in
\cat{F}\cap \cat{W}$ and cokernel $C\in \cat{C}$, and taking the pushout
instead of the pullback.  

Now, for an arbitrary map $f$, we write it as the composite 
\[
A\xrightarrow{i_{1}} A\oplus B \xrightarrow{f+1_{B}} B
\]
of a monomorphism followed by an epimorphism.  Write $f+1_{B}=q'j'$,
where $q'$ is an acyclic fibration and $j'$ is a cofibration.  Then write
$j'i_{1}=q''j$, where $q''$ is an acyclic fibration and $j$ is a
cofibration. Take $q=q''q'$ to complete the proof.  
\end{proof}

\begin{proposition}\label{prop-weak-composition}
Weak equivalences as defined above are closed under compositions.  
\end{proposition}

\begin{proof}
It suffices to check that a composition of the form $ip$, where $p$ is an
acyclic fibration and $i$ is an acyclic cofibration, can be written
$ip=qj$, where $q$ is an acyclic fibration and $j$ is an acyclic
cofibration.  By the preceding proposition, we can write $ip=qj$, where
$q$ is an acyclic fibration and $j$ is a cofibration.  This gives us the
diagram below 
\[
\begin{CD}
0 @>>> X @>j>> W @>>> \cok j @>>> 0 \\
@.  @VpVV @VqVV @VVrV \\
0 @>>> Y @>>i> Z @>>> \cok i @>>> 0
\end{CD}
\]
which leads to the short exact sequence 
\[
0 \xrightarrow{} \ker p \xrightarrow{} \ker q \xrightarrow{} \ker r
\xrightarrow{} 0.
\]
Since $p$ and $q$ are acyclic fibrations, $\ker p$ and $\ker q$ are in
$\cat{W}$.  Since $\cat{W}$ is thick, $\ker r\in \cat{W}$.  But $\cok
i\in \cat{W}$ since $i$ is an acyclic cofibration.  We conclude that
$\cok j\in \cat{W}$ since $\cat{W}$ is thick, and hence $j$ is an acyclic
cofibration as required.  
\end{proof}

\section{Cofibrant generation}

So now we have this correspondence between abelian model categories and
compatible pairs of complete cotorsion pairs.  We should then ask:
given an important property of model categories, how is that property
reflected in the compatible complete cotorsion pairs?  

For example, when is our abelian model structure cofibrantly generated?
Recall that a model structure is \textbf{cofibrantly generated}
(see~\cite[Section~3.1]{goerss-model}) when there is a set $I$ of
cofibrations and a set $J$ of acyclic cofibrations such that $p$ is an
acyclic fibration if and only if it has the \rlp $I$, and $p$ is a
fibration if and only if it has the \rlp $J$.  (There is also an
additional smallness condition that we omit, because it is automatically
satisfied in any standard algebraic category; it is only topologies that
make this one hard).  The key thing here is that we do not need the
entire proper class of cofibrations to detect the acyclic fibrations,
but just the \textbf{set} $I$.  Cofibrantly generated model categories
are much easier to work with than general model categories; for one
thing, the cokernels of the generating cofibrations play a somewhat
similar role as the spheres do in the algebraic topology of topological
spaces.  

The translation between abelian model structures and cotorsion pairs
basically takes a cofibration to its cokernel, so we define a cotorsion
pair $(\cat{D},\cat{E})$ to be \textbf{cogenerated by a set} when there
is a subset $\cat{D}'$ of the class $\cat{D}$ such that $E\in \cat{E}$
if and only if $\Ext^{1}(D,E)=0$ for all $D\in \cat{D}'$.  This
definition was actually made by Eklof and Trlifaj~\cite{eklof-trlifaj}
without knowing anything about model categories.

For example, the (projective, everything) cotorsion pair is cogenerated
by $0$ in any abelian category, and the (everything, injective)
cotorsion pair in the category of left $R$-modules is cogenerated by the
set of all $R/\ideal{a}$, where $\ideal{a}$ is a left ideal of $R$.
(This is Baer's criterion for injectivity).

Then the following lemma is not difficult. 

\begin{lemma}\label{lem-cofibrant}
If an abelian model category is cofibrantly generated, then the
corresponding complete cotorsion pairs $(\cat{C},\cat{F}\cap
\cat{W})$ and $(\cat{C}\cap \cat{W},\cat{F})$ are each cogenerated by a
set.
\end{lemma}

We would like the converse to be true as well.  In fact, we want more
than that.  Recall that the point of a model category being cofibrantly
generated is then Quillen's small object
argument~\cite[Theorem~3.5]{goerss-model} gives an automatic proof of
the factorization axioms that also proves the naturality of these
factorizations.  So we want to start with two compatible cotorsion
pairs, not necessarily complete, but cogenerated by a set, and argue
that the cotorsion pairs are automatically complete, and hence we get an
abelian model structure.  In fact, this seems to be true in practice,
but the simplest theorem along these lines requires a strong hypothesis.

\begin{proposition}\label{prop-cofibrant}
If $\cat{A}$ is a Grothendieck category with enough projectives, then
every cotorsion pair cogenerated by a set is complete.  Furthermore,
given two compatible cotorsion pairs each cogenerated by a set,
the corresponding abelian model structure is cofibrantly generated.  
\end{proposition}

The proof of this proposition can be found in~\cite{hovey-cotorsion}.  I
think of the Grothendieck hypothesis as the best hypothesis on an
abelian category.  It is general enough to include categories that occur
frequently in algebraic topology (sheaves and comodules, for example),
but strong enough to ensure good properties.  An abelian category is
Grothendieck when it is cocomplete and has a generator, and when filtered
colimits are exact.

The reason for having projectives is so that, given one of your
cogenerators $C$, you have a good choice for a monomorphism whose
cokernel is $C$.  Usually, even when you do not have enough projectives,
you actually do have such a good choice anyway, but it is more
complicated to make this into a theorem.  

\section{Monoidal structure}\label{sec-monoidal}

One of the most important properties a model structure can have is
compatibility with a tensor product.  This is particularly important in
the algebraic situation.  For example, given any Grothendieck category
$\cat{A}$, there is an injective model structure on unbounded chain
complexes over $\cat{A}$.  The cofibrations are the monomorphisms, the
weak equivalences are the homology isomorphisms, and the fibrations are
the epimorphisms with DG-injective kernel.  (DG-injective means each
entry is injective, and every map from an exact complex into it is chain
homotopic to $0$).  The homotopy category of the injective model
structure is the derived category of $\cat{A}$, and so the injective
model structure is the foundation for homological algebra of the $\Ext$
sort in any Grothendieck category.

But as a practical matter, one almost always has a tensor product
around; the tensor product of modules, or sheaves, or comodules.  And
injective resolutions are almost never compatible with the tensor
product, which means that one cannot use the injective model structure
to produce a derived tensor product on the derived category of
$\cat{A}$.  

In general, we have the following definition. 

\begin{definition}\label{defn-monoidal}
A model structure on a symmetric monoidal category $\cat{A}$ is called
\textbf{monoidal} whenever the following conditions hold: \begin{enumerate}
\item Given cofibrations $i\mathcolon A\xrightarrow{}B$ and $j\mathcolon
C\xrightarrow{}D$, the induced map 
\[
i\boxprod j\mathcolon (A\otimes D)\amalg_{A\otimes C}(B\otimes C)
\xrightarrow{} B\otimes D
\]
is a cofibration, which, in addition, is an acyclic cofibration if either
$i$ or $j$ is acyclic.  
\item An annoying condition that only arises when the unit of the tensor
product is not cofibrant (see~\cite[Definition~4.2.6]{hovey-model}).
\end{enumerate}
\end{definition}

The main point of a monoidal model category is that it gives the tensor
product on $\cat{A}$ homotopy-theoretic meaning.  Thus the homotopy
category of a monoidal model category $\cat{A}$ will itself be a
symmetric monoidal category, and one can usually also construct model
categories (and thus homotopy categories) of monoids in $\cat{A}$ and of
modules over a given monoid in $\cat{A}$.

The definition of a monoidal model category was not really even
formulated precisely until the late 1990's, although it is based on
Quillen's definition of a simplicial model category dating to the 1960's
(see~\cite[Definition~4.11]{goerss-model}).  Looking back on it,
however, one can say that one of the biggest problems in algebraic
topology was the failure to find a monoidal model category whose
homotopy category is the usual stable homotopy category.  This problem
was solved in the 1990's by Elmendorf, Kriz, Mandell, and
May~\cite{elmendorf-kriz-mandell-may} and
Smith~\cite{hovey-shipley-smith}.  

The projective model structure on chain complexes of $R$-modules is
monoidal, but, as mentioned above, the injective model structure is
not.  

Here is what a monoidal abelian model structure looks like from the
cotorsion pair point of view.

\begin{theorem}\label{thm-monoidal}
Let $\cat{A}$ be an abelian model category,and suppose $\cat{A}$ is
closed symmetric monoidal.  Suppose the following conditions are satisfied\uc 
\begin{enumerate}
\item Every element of $\cat{C}$ is flat. 
\item If $X,Y\in \cat{C}$, then $X\otimes Y\in \cat{C}$. 
\item If $X,Y\in \cat{C}$ and one of them is in $\cat{W}$, then
$X\otimes Y\in \cat{W}$.  
\item The unit $S$ is in $\cat{C}$. 
\end{enumerate}
Then $\cat{A}$ is a monoidal model category.  
\end{theorem}

Again, the proof of this theorem can be found in~\cite{hovey-cotorsion}.
Here ``flat'' means what it usually does.  That is, $X$ is flat if the
functor $X\otimes (-)$, which is right exact since it is a left adjoint,
is actually exact.

\section{Standard examples}

Having done the work of relating abelian model categories to pairs of
complete cotorsion pairs, we now consider the standard examples of
model structures on abelian categories.  

Perhaps the simplest example of a model category is the category of
$R$-modules when $R$ is a quasi-Frobenius ring.  This means that
projective and injective modules coincide.  The standard example is the
group ring $R=k[G]$ of a finite group $G$ over a field $k$.  In this
case, we can take $\cat{C}=\cat{F}$ to be the entire category of
$R$-modules, and take $\cat{W}$ to be the class of projective
(=injective) modules, which is thick in this unusual case.  The two
complete cotorsion pairs are then (everything, projective=injective) and
(projective=injective, everything).  The homotopy category of this model
category is called the \textbf{stable category of $R$-modules} and is
the main object of study in modular representation theory (as practiced
by Benson, Carlson, and Rickard, for example).  Two modules $M$ and $N$
are isomorphic in the stable category if there are projective modules
$P$ and $Q$ with $M\oplus P\cong N\oplus Q$.  The map
$M\xrightarrow{}M\oplus P$ is a typical acyclic cofibration and the map
$N\oplus Q\xrightarrow{}N$ is a typical acyclic fibration, but it is a
bit difficult to say exactly what a stable equivalence is other than a
composite of an acyclic cofibration followed by an acyclic fibration.  

Now suppose $\cat{A}$ is a Grothendieck category.  As mentioned above,
there is an injective model structure on $\Ch{\cat{A}}$, the category of
unbounded chain complexes on $\cat{A}$.  Here $\cat{C}$ is everything,
$\cat{W}$ is the exact complexes, and $\cat{F}$ consists of the
DG-injective complexes (defined at the beginning of
Section~\ref{sec-monoidal}).  A complex is DG-injective and exact if and
only if it is actually injective, so the cotorsion pair
$(\cat{C},\cat{F}\cap \cat{W})$ is (everything, injective).  The
cotorsion pair $(\cat{C}\cap \cat{W},\cat{F})$ is (exact,
DG-injective).  The homotopy category of the injective model structure
is the derived category of $\cat{A}$.  

The dual thing works for $\Ch{R}$, for $R$ a ring.  That is, we take
$\cat{C}$ to be the DG-projective complexes (defined after
Definition~\ref{defn-abelian}), $\cat{F}$ to be everything, and
$\cat{W}$ to be exact complexes.  Again, something that is both
DG-projective and exact is actually projective.  The homotopy category
of the projective model structure is the derived category of $R$, just
like the injective model structure.  However, the projective model
structure is monoidal when $R$ is commutative and so gives us more
structure (a derived tensor product and derived $\Hom $) on the derived
category than was apparent from the injective model structure.  

\section{Gorenstein rings}

Here is a new example of an abelian model category
from~\cite{hovey-cotorsion}.  The idea here is that we would like to do
modular representation theory over the integers instead of over fields.
So we want to study $\Z [G]$, when $G$ is a finite group.  This is no
longer a quasi-Frobenius ring; projectives and injectives do not coincide.
However, it does have some exceptionally nice properties: it is left and
right Noetherian, and $\Z [G]$, while not self-injective, does have
finite injective dimension as either a left or right module over
itself.  (This was first noticed by Eilenberg and
Nakayama~\cite{eilenberg-nakayama}).  Such a ring is called a
\textbf{Gorenstein ring}, or an \textbf{Iwanaga-Gorenstein ring}.  It is
a reasonable generalization of the usual notion of a commutative Gorenstein
ring.  

The salient fact about Gorenstein rings is that in a Gorenstein ring,
the modules of finite projective dimension and the modules of finite
injective dimension coincide (and the maximum injective or projective
dimension is the injective dimension of $R$).  This is due to
Iwanaga~\cite{iwanaga-1}.  It is easy to prove from this that these
modules form a thick subcategory.  We then define $\cat{W}$ to be this
class of modules with finite projective dimension, in analogy to the
quasi-Frobenius case.

But now the analogy breaks down a little, as we cannot expect to get a
model structure in which every module is both cofibrant and fibrant.  If
we want every module to be fibrant, then we take $\cat{C}$ to be the
class of \textbf{Gorenstein projective} modules; these are, of course,
modules $P$ for which $\Ext^{1}(P,W)=0$ for all $W$ of finite projective
dimension.  We should point out that Gorenstein projective modules are
still interesting over more general rings, but then the definition is
necessarily more complex~\cite[Definition~10.2.1]{enochs-jenda-book}.
Let $d$ be the injective dimension of $R$. Then a typical Gorenstein
projective is a $d$th syzygy of an arbitrary module.  That is, if we
take a module $M$ and take a partial projective resolution
\[
0 \xrightarrow{} K \xrightarrow{} P_{d-1} \xrightarrow{}\dotsb
\xrightarrow{}P_{0}\xrightarrow{}M \xrightarrow{} 0
\]
where the $P_{i}$ are projective, then $K$ is Gorenstein projective.
These modules have been studied before; when they are finitely
generated, they are called \textbf{maximal Cohen-Macaulay modules}.  

There is a dual notion of a \textbf{Gorenstein injective} module.  Here
$I$ is Gorenstein injective if and only if $\Ext^{1}(W,I)=0$ for all $W$
of finite projective dimension.  Again, this is only an appropriate
definition when the ring itself is Gorenstein; for the general case
see~\cite[Definition~10.1.1]{enochs-jenda-book}.  A typical Gorenstein
injective module is a $d$th cosyzygy of an arbitrary module.

We then get, after some work of course, two model structures on the
category of $R$-modules when $R$ is Gorenstein.  Both model categories
have the same class of acyclic objects $\cat{W}$, the modules of finite
projective dimension.  In the projective model structure, everything is
fibrant, and $M$ is cofibrant if and only if it is Gorenstein
projective.  In the injective model structure, everything is cofibrant,
and $M$ is fibrant if and only if it is Gorenstein injective.

The resulting homotopy category (which is the same for both model
structures) has every right to be called the \textbf{stable category of
$R$-modules}.  It is a triangulated category, and when $R=K[G]$ and $K$
is a principal ideal domain, it has a good closed symmetric monoidal
structure (given by tensoring over $K$).  It is the natural home for
representation theory of $G$ over $K$.  

As far as I know, not very much is known about this stable module
category.  There are many results about the stable module category of
$k[G]$ when $k$ is a field, such as a classification of the thick
subcategories of small objects (=finitely generated modules) when $G$ is
a $p$-group~\cite{benson-carlson-rickard-thick}.  It would be good to
know how much different the classification over $K[G]$ is.

\section{Gillespie's work}

The results in this section are due to my student, Jim Gillespie, and
come from~\cite{gillespie-modules},~\cite{gillespie-sheaves}, and
personal communications.  

\subsection{The general approach}

Gillespie looks at the general question of the relationship between a
cotorsion pair on a Grothendieck category $\cat{A}$ and the
homological algebra of $\cat{A}$.  That is, given a single cotorsion
pair $(\cat{D},\cat{E})$, can we induce a model structure on
$\Ch{\cat{A}}$ from this cotorsion pair on $\cat{A}$?  We know two
cases of this already: the (projective, everything) cotorsion pair on
$\cat{A}$ corresponds to the projective model structure on
$\Ch{\cat{A}}$, when it exists, and the (everything, injective)
cotorsion pair on $\cat{A}$ corresponds to the injective model
structure on $\Ch{\cat{A}}$.  

Recall how this works for the (projective, everything) model structure.
The two cotorsion pairs on $\Ch{\cat{A}}$ in this case are
(projective, everything) and (DG-projective, exact).  There is of
course a categorical definition of projective in $\Ch{\cat{A}}$, but
that will be of no help for a more general cotorsion pair.  Instead,
note that a complex $X$ is projective if and only if $X$ is exact and
$Z_{n}X$ is projective for all $n$.  This suggests the following
definition. 

\begin{definition}\label{defn-gillespie}
Suppose $\cat{D}$ is a class of objects in a bicomplete
abelian category $\cat{A}$.  Define $\widetilde{\cat{D}}$ to be the
class of objects $X$ in $\Ch{\cat{A}}$ such that $X$ is exact and
$Z_{n}X\in \cat{D}$ for all $n$. 
\end{definition}

So if $\cat{D}$ is projectives, we recover the notion of a projective
complex.  If $\cat{D}$ is everything, we recover the notion of an exact
complex.  

We still have to recover the notion of DG-projective.  Recall that $X$ is
DG-projective if each $X_{n}$ is projective and any map from $X$ to an
exact complex is chain homotopic to $0$.  This suggests the following
definition. 

\begin{definition}\label{defn-dg-gillespie}
Suppoe $(\cat{D},\cat{E})$ is a cotorsion pair in a bicomplete abelian
category $\cat{A}$.  Define $\dg \widetilde{\cat{D}}$ to be the class of
all $X$ in $\Ch{\cat{A}}$ such that $X_{n}\in \cat{D}$ for all $n$ and
every map from $X$ to a complex in $\widetilde{\cat{E}}$ is chain
homotopic to $0$.  Similarly, define $\dg \widetilde{\cat{E}}$ to be the
class of all $X\in \Ch{\cat{A}}$ such that $X_{n}\in \cat{E}$ for all
$n$ and every map from a complex in $\widetilde{\cat{D}}$ to $X$ is
chain homotopic to $0$.  
\end{definition}

So if $(\cat{D},\cat{E})$ is (projectives, everything), then $\dg
\widetilde{\cat{D}}$ is the class of DG-projectives and $\dg
\widetilde{\cat{E}}$ is everything.  Similarly, if $(\cat{D},\cat{E})$
is (everything, injectives), then $\widetilde{\cat{D}}$ is the class of
exact complexes, $\widetilde{\cat{E}}$ is the class of injective
complexes, $\dg \widetilde{\cat{D}}$ is everything, and $\dg
\widetilde{\cat{E}}$ is the class of DG-injective complexes. 

Now, the goal of Gillespie's work is to prove a metatheorem of the
following sort:

\begin{theorem}\label{thm-meta}
If $(\cat{D},\cat{E})$ is a nice enough cotorsion pair on a Grothendieck
abelian category $\cat{A}$, then there is an induced abelian model
structure on $\Ch{\cat{A}}$, where $\cat{C}=\dg \widetilde{\cat{D}}$,
$\cat{F}=\dg \widetilde{\cat{E}}$, and $\cat{W}$ is the class of exact
complexes.  
\end{theorem}

Of course, he also wants to give nontrivial examples of this theorem.  

Note that, because $\cat{W}$ is the category of exact complexes, the
homotopy category of any the model structures produced by
Theorem~\ref{thm-meta} is the usual derived category of $\cat{A}$.  So
this theorem is not producing new homotopy categories; instead, it is
producing new ways to understand the derived category.  This is
important if one wants the derived category to have some good properties
not accessible through the usual injective model structure.

\subsection{Making the theorem concrete}

We now need to specify precisely what it means for a cotorsion pair
$(\cat{D},\cat{E})$ to be ``nice enough'' in Theorem~\ref{thm-meta}.  

The first thing to verify is that $(\widetilde{\cat{D}},\dg
\widetilde{\cat{E}})$ and $(\dg
\widetilde{\cat{D}},\widetilde{\cat{E}})$ are indeed cotorsion
pairs.  This is simple enough that we can do it here, for
$(\widetilde{\cat{D}},\dg \widetilde{\cat{E}})$.  

We first show that $\Ext^{1}(Y,X)=0$ for $Y\in \widetilde{\cat{D}}$ and
$X\in \dg \widetilde{\cat{E}}$.  So suppose we have a short exact
sequence of complexes
\[
0 \xrightarrow{} X \xrightarrow{} W \xrightarrow{} Y \xrightarrow{} 0
\]
with $X\in \dg \widetilde{\cat{E}}$ and $Y\in \widetilde{\cat{D}}$.
Then each $X_{n}$ is in $\cat{E}$ and each $Y_{n}$ is in $\cat{D}$
(because $Z_{n}Y\in \cat{D}$ for all $n$ and $Y$ is exact, so $Y_{n}$ is
an extension of $Z_{n}Y$ and $Z_{n-1}Y$).  Therefore, our short exact
sequence of complexes is dimensionwise split, so $W_{n}\cong X_{n}\oplus
Y_{n}$.  In terms of this decomposition, the differential on $W$ is
$d=(d_{X},\tau +d_{Y})$, where $\tau \mathcolon
Y_{n}\xrightarrow{}X_{n-1}$.  Because $d^{2}=0$, we see that
$\tau\mathcolon Y\xrightarrow{}\Sigma X$ is a chain map.  By
hypothesis, this chain map is chain homotopic to $0$.  The chain
homotopy can then be used to define a splitting of our sequence by a
chain map.  

Now suppose $\Ext^{1}(Y,X)=0$ for all $X\in \dg \widetilde{\cat{E}}$.
We want to show that $Y\in \widetilde{\cat{D}}$.  The first thing to
point out is that 
\[
\Ext^{1}(Y, D^{n+1}A) \cong \Ext^{1}(Y_{n},A).  
\]
(To see this, just draw what an extension of complexes looks like).  It
follows easily from this that $Y_{n}\in \cat{D}$ for all $n$, since
$D^{n}A\in \dg \widetilde{\cat{E}}$ whenever $A\in \cat{E}$.

Given this, an element of $\Ext^{1}(Y, S^{n-1}A)$ is determined by a map
\[
Y_{n}/B_{n}Y\xrightarrow{}A
\]
(this is the same as a chain map $Y\xrightarrow{}S^{n}A$).  However, two
maps determine the same extension if they are chain homotopic as chain
maps $Y\xrightarrow{}S^{n}A$.  Said another way, $\Ext^{1}(Y, S^{n-1}A)$
is the quotient 
\[
\Hom (Y_{n}/B_{n}Y,A)/\Hom (Y_{n-1},A).
\]
If this quotient is to be $0$ for all $A\in \cat{E}$, we can in
particular take $A$ to be an injective object containing $Y_{n}/B_{n}Y$
to see that $Y$ is exact.  But then $\Ext^{1}(Y,S^{n-1}A)$ is isomorphic
to $\Ext^{1}(Z_{n-1}Y,A)$, from which we see that $Z_{n-1}Y\in \cat{D}$
for all $n$.  Thus $Y\in \widetilde{\cat{D}}$ as required.

A similar, but simpler, argument shows that if $\Ext^{1}(Y,X)=0$ for all
$Y\in \widetilde{\cat{D}}$, then $X\in \dg \widetilde{\cat{E}}$.  For
this, one uses the isomorphism 
\[
\Ext^{1}(D^{n}A, X)\cong \Ext^{1}(A,X_{n})
\]
to see that $X_{n}\in \cat{E}$ for all $n$.  It then follows that any
element in $\Ext^{1}(Y,X)$ is dimensionwise split for $Y\in
\widetilde{\cat{D}}$, so $\Ext^{1}(Y,X)$ is isomorphic to chain
homotopy classes of chain maps from $Y$ to $\Sigma X$.  Since
$\Ext^{1}(Y,X)=0$, we see that $X\in \dg \widetilde{\cat{E}}$.  

Now, if we worked with $(\dg \widetilde{\cat{D}},\widetilde{\cat{E}})$
instead, we would have run into a problem.  In the above argument, there
was a point where we embedded $Y_{n}/B_{n}Y$ into an element of
$\cat{E}$, which we can do by taking an injective.  The dual will cause
us trouble because we do not want to assume there are enough projectives
in $\cat{A}$.  So instead we assume there are enough $\cat{D}$-objects
in $\cat{A}$, in the sense that everything in $\cat{A}$ is a quotient of
something in $\cat{D}$.  This would be automatic if $(\cat{D},\cat{E})$
were a complete cotorsion pair.  

So we get the following proposition of Gillespie.  

\begin{proposition}
If $(\cat{D},\cat{E})$ is a cotorsion pair on a Grothendieck category
$\cat{A}$ that has enough $\cat{D}$-objects, then
$(\widetilde{\cat{D}},\dg \widetilde{\cat{E}})$ and $(\dg
\widetilde{\cat{D}},\widetilde{\cat{E}})$ are cotorsion pairs on
$\Ch{\cat{A}}$.
\end{proposition}

We now want to know whether these cotorsion pairs are compatible with
the class $\cat{W}$ of exact complexes.  That is, we want to know that 
\[
\dg \widetilde{\cat{D}}\cap \cat{W}=\widetilde{\cat{D}} \text{ and } \dg
\widetilde{\cat{E}}\cap \cat{W} =\widetilde{\cat{E}}. 
\]
It is fairly straightforward to show the inclusions 
\[
\widetilde{\cat{D}}\subseteq \dg \widetilde{\cat{D}}\cap \cat{W} \text{
and } \widetilde{\cat{E}} \subseteq \dg \widetilde{\cat{E}}\cap \cat{W}.
\]
One shows that any map from something in $\widetilde{\cat{D}}$ to
something in $\widetilde{\cat{E}}$ is chain homotopic to $0$.  The idea
for the converse is as follows.  Given $X\in \dg \widetilde{\cat{D}}\cap
\cat{W}$, we want to show $\Ext^{1}(Z_{n}X,A)=0$ for all $A\in
\cat{E}$.  Since we have the short exact sequence 
\[
0 \xrightarrow{} Z_{n+1}X \xrightarrow{} X_{n+1} \xrightarrow{}
Z_{n}X\xrightarrow{} 0
\]
and $X_{n+1}\in \cat{D}$, it suffices to show that any map
$Z_{n+1}X\xrightarrow{}A$ extends to a map $X_{n+1}\xrightarrow{}A$.
Take an augmented injective resolution $I_{*}$ of $A$ (so $I_{0}=A$ and
$I_{-1}$ is an injective object containing $A$).  With any justice, this
should be a complex in $\widetilde{\cat{E}}$, since $A$ was in $\cat{E}$
to start with.  Then a map $Z_{n+1}X\xrightarrow{}A$ induces a map of
complexes $\Sigma^{-n-2}X\xrightarrow{}I_{*}$ using injectivity.  This
chain map is chain homotopic to $0$, and the chain homotopy gives us an
extension $X_{n+1}\xrightarrow{}A$.  

This argument depended on $I_{*}$ actually being in
$\widetilde{\cat{E}}$.  This is \textbf{NOT} automatic, however.
Consider the following three conditions on a cotorsion pair $(\cat{D},
\cat{E})$.  

\begin{enumerate}
\item $\Ext^{i}(D,E)=0$ for all $D\in \cat{D}$, $E\in \cat{E}$, and
$i>0$.  
\item $\cat{D}$ is closed under kernels of epimorphisms.  
\item $\cat{E}$ is closed under cokernels of monomorphisms.  
\end{enumerate}

One can check easily that the first condition above implies the second
and third; when our cotorsion pair satisfies this first condition, we
call it a \textbf{hereditary} cotorsion pair.  The second condition is
equivalent to the first when our category has enough projectives, and
the third condition is equivalent to the first when our category has
enough injectives.  

Most cotorsion pairs that arise naturally are hereditary, though it can
sometimes be hard to prove that a cotorsion pair is hereditary if there
are not enough projectives in the category.  

Then we have the following proposition, again due to Gillespie. 

\begin{proposition}\label{prop-gillespie-hereditary}
Suppose $(\cat{D},\cat{E})$ is a hereditary cotorsion pair in a
Grothendieck category $\cat{A}$ with enough $\cat{D}$-objects.  Then
$\dg \widetilde{\cat{D}}\cap \cat{W}=\widetilde{\cat{D}}$.  If, in
addition, $\cat{A}$ has enough projectives, then $\dg
\widetilde{\cat{E}}\cap \cat{W}=\widetilde{\cat{E}}$.  
\end{proposition}

As a practical matter, though, most of the categories we are interested
in do not have enough projectives.  Gillespie and Ed Enochs get around
this with a subtle transfinite induction argument that gives the
following proposition.

\begin{proposition}\label{prop-gillespie-transfinite}
Suppose $(\cat{D},\cat{E})$ is a cotorsion pair that is cogenerated by
a set on a Grothendieck category $\cat{A}$ with enough
$\cat{D}$-objects.  Then $(\dg \widetilde{\cat{D}},\widetilde{\cat{E}})$
has enough injectives.  
\end{proposition}

One can look on this proposition as a variant of the small object
argument, but it is much more complicated to prove.  Also, it does not
seem to work for $(\widetilde{\cat{D}},\dg \widetilde{\cat{E}})$.   

From this, then, we get the following proposition of Gillespie.  

\begin{proposition}\label{prop-gillespie-almost}
Suppose $(\cat{D},\cat{E})$ is a hereditary torsion theory cogenerated
by a set on a Grothendieck category $\cat{A}$ with enough
$\cat{D}$-objects.  Then 
\[
\dg \widetilde{\cat{D}}\cap \cat{W}=\widetilde{\cat{D}}, \dg
\widetilde{\cat{E}}\cap \cat{W}=\widetilde{\cat{E}}, 
\]
and $(\dg \widetilde{\cat{D}},\widetilde{\cat{E}})$ is complete.  
\end{proposition}

The proof is not hard.  Suppose $X\in \dg \widetilde{\cat{E}}\cap
\cat{W}$.  We have a short exact sequence 
\[
0 \xrightarrow{} X \xrightarrow{} W \xrightarrow{}Y \xrightarrow{} 0
\]
with $W\in \widetilde{\cat{E}}$ and $Y\in \dg \widetilde{\cat{D}}$.  But
then $X$ and $W$ are exact, so $Y$ is too.  Thus $Y\in
\widetilde{\cat{D}}$.  This means the sequence splits, so $X$ is a
summand in $W$.  But then $X\in \widetilde{\cat{E}}$.  

Given that $(\dg \widetilde{\cat{D}},\widetilde{\cat{E}})$ has enough
injectives, we can use a pushout trick to show it has enough projectives
as well, using the fact that there are enough $\cat{D}$-objects.  That
is, you first show that $\Ch{\cat{A}}$ has enough
$\widetilde{\cat{D}}$-objects.  Then, given $X$, you take a surjection
$A\xrightarrow{}X$ with kernel $K$, where $A\in \widetilde{\cat{D}}$.
Then you embed $K$ in an element of $\widetilde{\cat{E}}$ with cokernel
in $\dg \widetilde{\cat{D}}$, and you take the pushout.  

So, to complete Gillespie's program, we must ensure that
$(\widetilde{\cat{D}},\dg \widetilde{\cat{E}})$ is complete.  In fact,
the pushout trick above shows that we only need to be sure it has enough
injectives.  This appears to be the heart of the matter.

One always wants to use some version of the small object argument of
Quillen.  But it just seems to be harder than it is for model
categories, and being cogenerated by a set does not seem to be enough.
So Gillespie, following Enochs and Lop\'{e}z-Ramos~\cite{enochs-lopez},
strengthens the definition a bit.

\begin{definition}\label{defn-Kaplansky}
A class $\cat{D}$ is a \textbf{Kaplansky class} if there is
some cardinal $\kappa$ such that, for every $\kappa$-generated subobject
$T$ of an object $D\in \cat{D}$, there is a $\kappa$-presentable object
$S\in \cat{D}$ such that $T\subseteq S\subseteq D$ and $D/S\in
\cat{D}$.  
\end{definition}

In an arbitrary category, an object $A$ is $\kappa$-generated if $\Hom
(A,-)$ commutes with $\lambda$-fold coproducts for all $\kappa$-filtered
ordinals $\lambda$ (any regular cardinal larger than $\kappa $ is
$\kappa$-filtered).  On the other hand, $A$ is $\kappa$-presentable if
$\Hom (A,-)$ commutes with all $\kappa$-filtered colimits.  The easiest
case is when $\kappa =\omega$, when we do recover the usual definition
of finitely generated and finitely presentable, only without reference
to a specific generator of the category.  

This is a strange definition at first.  It is motivated by flat
modules, where it asserts that, given any small subset of a flat
module, there is a flat submodule that contains it and sits inside the
big module purely.  This was the key idea in the proof of the flat cover
conjecture by Bican, El Bashir, and
Enochs~\cite{bican-el-bashir-enochs}.  

We can now state the precise version of Gillespie's main theorem.  

\begin{theorem}\label{thm-gillespie}
Suppose $(\cat{D},\cat{E})$ is a hereditary cotorsion pair cogenerated
by a set such that $\cat{D}$ is a Kaplansky class on a Grothendieck
category $\cat{A}$ with enough $\cat{D}$-objects.  Then there is an
induced abelian model structure on $\Ch{\cat{A}}$, where $\cat{C}=\dg
\widetilde{\cat{D}}$, $\cat{F}=\dg \widetilde{\cat{E}}$, and $\cat{W}$
is the class of exact complexes.
\end{theorem}

\subsection{Sheaves and schemes}

The motivation for Gillespie's work was to better understand the derived
category of sheaves on a ringed space and the derived category of
quasi-coherent sheaves on a scheme.  Recall that a ringed space is a
topological space $S$ equipped with a sheaf of rings $\cat{O}$; that is,
a contravariant functor from open sets of $S$ to commutative rings that
is locally determined (the sheaf property).  A one-point ringed space is
of course a commutative ring.  The category of $\cat{O}$-modules is
therefore a generalization of the category of $R$-modules for a
commutative ring $R$; here an $\cat{O}$-module $M$ is a sheaf of abelian
groups over $S$ such that $M(U)$ is naturally a module over
$\cat{O}(U)$.  The category of $\cat{O}$-modules has a lot in common
with the category of $R$-modules; it is a closed symmetric monoidal
Grothendieck category.  There is a tensor product defined stalkwise in
the obvious way.

Therefore, we would expect $\Ch{\cat{O}}$ to be a symmetric monoidal
model category, so that the derived category of $\cat{O}$-modules
inherits a tensor product.  However, before Gillespie's work, I don't
believe this was known.  The injective model structure certainly exists
on $\Ch{\cat{O}}$, but it is not compatible with the tensor product, and
cannot be used to define a derived tensor product.  The projective model
structure only exists rarely, because generally there are not enough
projective $\cat{O}$-modules.  There are enough flats, though; in fact,
the flat sheaves $\cat{O}_{U}$ generate the category, where the stalks
of $\cat{O}_{U}$ agree with the stalks of $\cat{O}$ inside $U$ and are
$0$ outside $U$.  The author used these sheaves to construct a monoidal
model structure on $\Ch{\cat{O}}$ in~\cite{hovey-sheaves}, but only
under an annoying technical assumption on the ringed space, involving
the finiteness of sheaf cohomology.  

Gillespie's work allows one to use all the flat sheaves at once, rather
than just the ones one can explicitly write down.  That is, we start
with the (flat, cotorsion) cotorsion pair on $\cat{O}$-modules.  One
needs an argument involving the stalks to see that this cotorsion pair
is hereditary.  Using the approach to the flat cover conjecture
of~\cite{bican-el-bashir-enochs}, Gillespie shows that that the flat
$\cat{O}$-modules form a Kaplansky class.  (The proof involves purity in
an essential way).  Hence Theorem~\ref{thm-gillespie} gives us an
abelian model structure on $\Ch{\cat{O}}$, which Gillespie proves is
compatible with the tensor product.  Gillespie therefore gets a derived
tensor product and a derived $\Hom$ functor with all the usual
properties on the derived category of $\cat{O}$-modules.

In algebraic geometry, however, it is more common to use the category of
quasi-coherent sheaves on a scheme.  This is because if your ringed
space is $\text{Spec } R$, then a quasi-coherent sheaf is equivalent to
an $R$-module, whereas an arbitrary sheaf could be more complicated.
The word quasi-coherent can be thought of as meaning locally a quotient
of free sheaves.  We need a special argument to see that the (flat,
cotorsion) pair is hereditary, involving comparision to open affine
subschemes, that only works if the scheme is quasi-separated.  Now
Enochs, Estrada, Garc\'{a} Rozas, and
Oyonarte~\cite{enochs-estrada-garcia-oyonarte} have proved a result
equivalent to the fact that flat quasi-coherent sheaves form a Kaplansky
class (see also~\cite[Proposition~3.3]{enochs-estrada-flat}).  There is
an additional complication though; it is much less obvious that there
are enough flat quasi-coherent sheaves (the sheaves $\cat{O}_{U}$ are
not quasi-coherent).  This is known by algebraic geometers when the
scheme is quasi-compact and separated, and Gillespie and I suspect it
holds when the scheme is quasi-compact and quasi-separated (this seems
to be the hypothesis of choice in algebraic geometry anyway).  But in
any case, Gillespie's work then leads to an abelian monoidal model
structure of chain complexes of quasi-coherent sheaves over a
quasi-compact, separated scheme, and hence a derived tensor product and
a derived $\Hom$ functor on the derived category of the scheme.  This
derived tensor product is used frequently by algebraic geometers, and
this provides a simple reason for its existence and a simple proof that
it has all the properties one would expect.  (I believe the usual
approach is to patch together the derived tensor products of each affine
piece of the scheme).


\providecommand{\bysame}{\leavevmode\hbox to3em{\hrulefill}\thinspace}
\providecommand{\MR}{\relax\ifhmode\unskip\space\fi MR }
\providecommand{\MRhref}[2]{%
  \href{http://www.ams.org/mathscinet-getitem?mr=#1}{#2}
}
\providecommand{\href}[2]{#2}

\end{document}